\documentclass[12pt]{amsart}
\usepackage{amssymb,latexsym,eufrak,amsmath,amscd}
\usepackage[all]{xy}
\theoremstyle{plain}
\newtheorem{theorem}{Theorem}[section]
\newtheorem{proposition}[theorem]{Proposition}

\newtheorem{lemma}[theorem]{Lemma}

\newtheorem*{variantnn}{Variant}

\newtheorem*{theoremnna}{Theorem A}
\newtheorem*{theoremnnb}{Theorem B}

\theoremstyle{definition}
\newtheorem{definition}[theorem]{Definition}
\newtheorem{remark}[theorem]{Remark}
\newtheorem{example}[theorem]{Example}

\newcommand{\lra}{\longrightarrow}
\newcommand{\noi}{\noindent}
\newcommand{\PP}{\mathbf{P}}

\newcommand{\NN}{\mathbf{N}}

\newcommand{\CC}{\mathbf{C}}
\newcommand{\QQ}{\mathbf{Q}}

\newcommand{\Zero}[1]{\text{Zeroes}\big({#1}\big)}
\newcommand{\OO}{\mathcal{O}}

\newcommand{\II}{\mathcal{I}}

\newcommand{\fra}{\frak{a}}
\newcommand{\frb}{\frak{b}}
\newcommand{\fro}{\frak{o}}
\newcommand{\frmm}{\frak{m}}
\newcommand{\frq}{\frak{q}}
\newcommand{\frr}{\frak{r}}
\newcommand{\frp}{\frak{p}}

\newcommand{\HH}[3]{H^{{#1}} \big( {#2} , {#3} \big) }

\newcommand{\for}{ \ \ \text{ for } \ }
\newcommand{\fall}{ \ \ \text{ for all } \ }

\newcommand{\MI}[1]{\mathcal{J} ( {#1} ) }

\newcommand{\ord}{\text{ord}}

\newcommand{\pr}{\prime}

\newcommand{\Bl}{\text{Bl}}

\title{Uniform Bounds and Symbolic Powers on  Smooth
Varieties} 

\author{Lawrence Ein}
\address{Department of Mathematics \\University of
Illinois at Chicago \hfil\break\indent 
851 South Morgan Street (M/C 249)\\ Chicago, IL 60607-7045, USA}
\email{ein@math.uic.edu}
\thanks{Research of the first author partially
supported by NSF Grant DMS  99-70295}

\author{Robert Lazarsfeld} 
\address{Department of Mathematics
\\ University of Michigan \\ Ann Arbor, MI 48109, USA}
\email{rlaz@math.lsa.umich.edu}
\thanks{Research of the second author  partially
supported by   NSF Grant DMS
97-13149}

\author{Karen E. Smith}
\address{Department of Mathematics
\\ University of Michigan \\ Ann Arbor, MI 48109, USA}
\email{kesmith@math.lsa.umich.edu}
\thanks{Research of the  third author  partially
supported by  NSF Grant DMS 96-25308}

\begin{document}

\maketitle
%\setcounter{section}{-1}
%\tableofcontents

%\setlength{\parskip}{.1 in}

\section*{Introduction}

The purpose of this note is to show how one can use 
multiplier ideals to establish effective uniform bounds
on the multiplicative behavior of certain families of
ideal sheaves on a smooth algebraic variety. In
particular, we prove a quick but rather surprising
result concerning the symbolic powers of radical ideals
on   such a variety.

Let $X$ be a non-singular quasi-projective
variety defined over the complex numbers $\CC$, and
let $Z \subseteq X$ be a reduced subscheme of
$X$.\footnote{All of our results are local in nature,
so there is no loss in taking $X$ to be an affine
variety. In this case one can work with the
coordinate ring $\CC[X]$ of $X$ in place of its
structure sheaf $\OO_X$.} 
Denote by  
\[ \frq \ = \ \II_Z  \ \subseteq \ \OO_X \]
the ideal sheaf of $Z$, so that $\frq$ is a sheaf of
radical ideals. We shall be concerned with the
symbolic powers $\frq^{(m)}$ of $\frq$. 
According
to a well-known theorem of Zariski and Nagata
(see \cite{Eisenbud}, Chapter 3.9)  $\frq^{(m)}$ can
be described as  the sheaf of all function germs
vanishing to order
$\ge m$ at a general point of each irreducible
component of
$Z$ (or equivalently at every point of
$Z$): 
\[ \frq^{(m)}  \ = \ \big \{   
\ f \in \OO_X \mid
\ord_x(f)
\ge m \ \text{for all $x \in Z$} \ \big \}. \]

It is evident that  $\frq^m \subseteq \frq^{(m)}$, but
in general of course the inclusion is strict. However
Swanson
\cite{Swanson}    established (in a much less
restrictive setting\footnote{Swanson's theorem holds in
particular on any normal variety over a field of any
characteristic.}) that there exists an integer
$k = k(Z)$ depending on $Z$ such that 
\[ \frq^{(km)} \ \subseteq  \ \frq^m  \fall m \in 
\NN. \]
On geometric grounds this already seems rather
striking since membership in the symbolic power on the
left is tested at   general smooth points of $Z$,
whereas the actual power on the right reflects also its
singular points. So one's first guess might be that the
worse the singularities of $Z$, the larger one will
have to take the coefficient
$k(Z)$ to be. Surprisingly enough  
 this is not the case, and in fact  
one has a uniform statement depending only on the
codimension of $Z$:
\begin{theoremnna}
Assume that every component of $Z$ has codimension $\le
e$ in $X$. Then
\[ \frq^{(me)} \ \subseteq \ \frq^m \ \fall m \in
\NN.\]
\end{theoremnna}
\noi In particular, if $\dim X = n$  then $\frq^{(mn)}
\subseteq \frq^m  $  for every radical
ideal
$\frq
\subseteq \OO_X$  and every natural number $m \ge 1$. 
One can see the Theorem as providing further
confirmation of  Huneke's philosophy
\cite{Huneke} that there are unexpected uniform bounds
lurking in commutative algebra.

Theorem A is a very simple application of the theory
of multiplier ideals.  In commutative algebra these
were introduced   by Lipman
\cite{Lipman} in connection with the Brian\c con-Skoda
theorem.\footnote{Lipman called them ``adjoint
ideals", but ``multiplier ideal" has become  standard
in higher dimensional geometry. The name derives from
their analytic construction, where
they arise as sheaves of multipliers (see
\cite{Demailly}).} More general
constructions, which we use here, have in the meantime
become extremely important in the study of higher
dimensional algebraic varieties (cf. \cite{Demailly},
\cite{Demailly.Bourbaki}, \cite{Ein}, \cite{Siu},
\cite{Kawamata},
\cite{MIAG}). In brief, we  consider   families
$
\fra_{\bullet} = \{ \fra_k
\}$ of ideals
$\fra_k \subseteq \OO_X$  --- such as the symbolic
powers
$\frq^{(\bullet)} = \{ \frq^{(k)} \} $ --- satisfying
the relations
\[ \fra_{\ell} \cdot \fra_{m} \ \subseteq \ 
\fra_{\ell + m} 
\fall \ell , m \ge 1. \]
 For each $\ell \ge 1$ we
associate to such a family
  an   \textit{asymptotic multiplier
ideal} 
$ \MI{ \Vert \fra_{\ell} \Vert}
\subseteq \OO_X $ which reflects the
asymptotic behavior of all the ideals
$\fra_{ p\ell }$ for $p \gg 0$. Using the subadditivity
theorem of
\cite{DEL},  we prove
\begin{theoremnnb}
If $\MI{\Vert \fra_{\ell}
\Vert} \subseteq \frb$ for some index $\ell$ and some
ideal $\frb$, then 
\[
\fra_{m\ell }
 \ \subseteq \  \frb^{
 m } \] for every $m \ge 1$.
\end{theoremnnb} 
\noi  In the case of symbolic
powers it is elementary to check that $\MI{\Vert
\frq^{(e)} \Vert} \subseteq
\frq$, so  Theorem A follows from the ``abstract"
Theorem B. As another application, we establish a
result (Theorem \ref{Izumi}) rendering effective and
extending in certain directions a theorem of Izumi 
\cite{Iz},
\cite{Hub.Swan} dealing with  ideals arising from  a
valuation.  

 We have been guided by the viewpoint  that
the families $\{ \fra_k \}$ share some of the
behavior of the linear series
$| k D |$ associated to multiples of a divisor 
$D$ on a projective variety,  and that one can try to
adapt geometric tools to the present setting. We hope
that these and other ideas from higher dimensional
complex geometry will find further algebraic
applications in the future. Going in the other
direction,   Hochster and Huneke \cite{HH} have used the
theory of tight closure to reprove and  generalize  
the    results of the present paper dealing with  
symbolic powers: in paticular they show that Theorem A
holds for any  regular local ring containing a
field, and they remove the hypothesis that $\frq$ be
radical (see
\S 3 for further
discussion).  This illustrates once again the close but
somewhat mysterious connections between tight closure
methods and the more geometric outlook appearing here.

Our exposition is organized into two sections.
In the first, we construct the multiplier ideals we
use and establish their basic properties. The 
applications are given in
\S 2.

We are grateful to Mel Hochster and Craig Huneke for
valuable discussions and encouragement, and to Jessica
Sidman for some Macaulay scripts related to Example
\ref{Pts.in.Plane}.  We also wish to record our debt to
the work of Irena Swanson and her collaborators,
through which we learned of many of the  questions
discussed here. 

\section{Graded Families and Multiplier Ideals}

In this section we construct the multiplier ideals we
require, and give  their basic
properties. Quick overviews of the general theory of
multiplier ideals appear in
\cite{Ein} and \cite{DEL}, \S 1, and a survey of some
of the applications in algebraic geometry is given in
\cite{Demailly.Bourbaki}.  The forthcoming book
\cite{PAG} will contain a detailed exposition, which in
the meantime can be found in the lecture notes
\cite{MIAG}. In particular, \cite{MIAG} contains full
proofs of all the facts about multiplier ideals quoted
in the following paragraphs.

Let  $X$ be a non-singular complex
quasi-projective variety, and 
$\fra \subseteq \OO_X$ a non-zero ideal sheaf on  $X$. A
\textit{log resolution} of $\fra$ is a projective
birational map $\mu : X^{\pr} \lra X$, with $X^{\pr}$
non-singular, such that $\fra \cdot \OO_{X^\pr} =
\OO_{X^\pr}(-F)$ for an effective Cartier divisor $F$ on
$X$ with the property that the sum of $F$ and the
exceptional divisor of $\mu$ has simple normal
crossing support. Such resolutions can be
construced by resolving the singularities of the
blow-up of $\fra$. We write $K_{X^\pr / X} = K_{X^\pr} -
\mu^* K_X$ for the relative canonical divisor of
$X^\pr$ over $X$.

  Given a rational number $c > 0$, the
\textit{multiplier ideal} associated to $c$ and $\fra$
is defined by fixing a log resolution as above, and
setting
\[ \MI{X , c \cdot \fra} \ = \ \MI{c \cdot \fra} \ = \
\mu_*
\OO_{X^\pr}\big(\, K_{X^\pr / X}  - [ \, c
\, F\, ]\, \big ).
\]
 Here $c F$ is viewed as an effective
$\QQ$-divisor, and its integer part $[ c \,F]$ is
defined by taking the integral part of the
coefficient of each of its components. The fact that 
$\MI{c
\cdot \fra}$ is indeed an ideal follows from the
observation that $\MI{c \cdot \fra} \subseteq \mu_*
\OO_{X^\pr}( K_{X^\pr / X}) = \OO_X$. An important
point is that this definition is independent of the log
resolution $\mu$. 

It follows immediately from the
definition that if $c \in \NN$, then $\MI{c \cdot \fra}
= \MI{\fra^c}$. This being so, we sometimes
prefer  to use ``exponential notation" $\MI{\fra^c}$
for the multiplier ideal $\MI{c \cdot \fra}$ for an
arbitrary rational number $c > 0$. Note that we are not
trying to attach any actual meaning to the expression $c
\cdot
\fra$ or $\fra ^ c$ when $c$ is non-integral.
Nonetheless, the possibility of being able to work with
rational coefficients is critical in applications. 

 As a variant,
given ideals $\fra , \frb \subseteq \OO_X$, and
rational numbers $c , d > 0$, we define $\MI{ (c \cdot
\fra) \cdot (d \cdot \frb)}$ (or $\MI{\fra^c \cdot
\frb^d}$ in exponential notation) by taking a common
log resolution
$\mu : X^\pr
\lra X$ of $\fra$ and $\frb$, with $ \fra \cdot
\OO_{X^\pr} =
\OO_{X^\pr}(-F_1),  \ \frb \cdot \OO_{X^\pr} =
\OO_{X^\pr}(-F_2)$, and putting $\MI{\fra^c \cdot
\frb^d} = \mu_* \OO_{X^\pr}\big(K_{X^\pr/X}-[cF_1 +
dF_2]\big)$. It is sometimes useful also to adopt the
convention that if $\fra = (0)$, then $\MI{c \cdot
\fra} = (0)$ for all $c > 0$. 

The most important local property of multiplier
ideals  is the
\textit{Restriction Theorem}, due in the
algebro-geometric setting to Esnault and Viehweg.
Specifically, let
$Y \subseteq X$ be a smooth subvariety, and let $\fra
\subseteq \OO_X$ be an ideal sheaf whose zeroes do not
contain $Y$. Then $\fra_Y = \fra \cdot \OO_Y$ is an
ideal sheaf on $Y$, and the result in question states
that one has an inclusion:
\begin{equation} \MI{\, Y \, , \, c \cdot \fra_Y} \
\subseteq
\
\MI{\, X \, , \, c \cdot \fra} \cdot \OO_Y 
\end{equation} of ideal sheaves on $Y$. This is
established by reducing to the case where
$Y$ has codimension one, and applying vanishing
theorems. The restriction theorem leads in turn to the
\textit{Subaddivity Theorem} of \cite{DEL}, which
states (in exponential notation) that given ideals $\fra
, \frb \subseteq \OO_X$  and  rational numbers
$c , d > 0$, one has the inclusion:
\begin{equation}  \label{Subadditivity} \MI{ \fra^c
\cdot
\frb^d}
\subseteq
\MI{\fra^c}
\cdot \MI{\frb^d}. \end{equation}
 To prove this, one first of all
applies the K\"unneth formula to check that 
\[ \MI{\, X \times X\, , \, (p_1^{-1} \fra)^c \cdot
(p_2^{-1}
\frb)^d \,} \ = \ p_1^{-1} \MI{\, X \, , \, \fra^c}
\cdot  p_2^{-1}
\MI{X, \frb^d \,},\]
where $p_1, p_2 : X \times X \lra X$ are the
projections, where we are somewhat abusively writing 
$f^{-1} \fro \subseteq \OO_V$ for inverse image $\fro
\cdot \OO_V$ of an ideal $\fro \subseteq \OO_W$ under a
morphism
$f : V
\lra W$. Then one restricts to the diagonal. Note for
later reference that in ``additive notation" 
(\ref{Subadditivity}) implies
\begin{equation} \label{Subadditivity.additive}
\MI{ cm \cdot \fra} \ \subseteq \ \MI{ c \cdot \fra}^m
\end{equation}
for every integer $m \ge 1$.

Some of the most interesting applications of multiplier
ideals (for instance \cite{Siu}, \cite{Kawamata})
depend on the fact that one can make asymptotic
constructions. A natural algebraic setting for these is
described in the following
\begin{definition} A \textit{graded family} or
\textit{graded system of ideals} 
$\fra_{\bullet} = \{ \fra_k \}$ is a collection of ideal
sheaves $\fra_k \subseteq \OO_X$ $(k \ge 1)$ satisfying
\begin{equation} \label{gsi.cond} 
\fra_k \cdot \fra_\ell \ \subseteq \ \fra_{k + \ell}
\fall  \ k, \ell \ge 1. 
\end{equation}
To avoid unnecessary complications, we 
assume also that $\fra_k \ne (0)$ for $k \gg 0$. 
\end{definition}

Note that if we
set $A_0 = \OO_X$,  then condition (\ref{gsi.cond})
is  equivalent to the statement that
$\oplus_{k \ge 0}
\ \fra_k$ is a graded $\OO_X$-algebra. The
asymptotic constructions that follow are particularly
useful in case this algebra is not finitely
generated (or at least not known to be so). 

\begin{example} \label{Examples.GSI} \begin{itemize}
\item[(i).] Let $(0) \ne \fra \subseteq \OO_X$ be a
fixed ideal, and take 
$\fra_k = \fra^k$ to be the $k^{\text{th}}$ power of
$\fra$. Then the $\{ \fra_k \}$ form a graded
family. One should view this as a trivial example.
\vskip 5pt
\item[(ii).] Let $D$ be a  divisor on a
projective variety $X$. When $\HH{0}{X}{\OO_X(kD)} \ne
0$ let
$\frb_k =
\frb \big( | kD | \big)$ be the base-ideal of the
complete linear series
$| kD |$, and put $\frb_k = (0)$ otherwise.
Then
$\frb_{\bullet} = \{ \frb_k \}$ forms a graded family
of ideals.

\vskip5pt
\item[(iii).] Let $(0) \ne \frq \subseteq \OO_X$ be a
radical ideal. Then the symbolic powers $\{\frq^{(k)}
\}$  form a graded family of ideals that we denote by
$\frq^{(\bullet)}$. 

\vskip 5pt
\item[(iv).]  Let $\nu : Y \lra X$ be a proper
birational map, and let $D$ be a non-zero 
effective Cartier divisor on $Y$. Then we get a graded
family of ideals
$\fro_{\bullet} = \{ \fro_k \}$ on $X$ by putting
$\fro_k =
\nu_* \OO_Y(-kD)$. Note that this includes the symbolic
powers $\frq^{(k)}$ in (iii) as a special case, as well
as the graded family of ideals associated to an
$\frmm$-valuation on $X$ in the sense of
\cite{Hub.Swan}. 

\vskip 5pt
\item[(v).] Let $p(t) = \sum_{i = 1}^{\infty}
\tfrac{1}{i!} t^i \in \CC[[t]]$ be the power series of
the function $e^t - 1$, and given $f \in \CC[x,y]$
define
\[ v(f)  \ = \ \text{ord}_t  \,  f\big(\, t\, , \, p(t)
\, \big).
\] This is a valuation on $\CC[x,y]$, and therefore the
ideals 
\[ \fro_k \  =_{\text{def}} \ \big \{ f \in \CC[x,y]
\mid v(f) \ge k \big \} \]
(which we may view as ideal sheaves on $X = \CC^2$) form
a graded family. Explicitly, 
\[\fro_k \ = \big( \ x^k  \ , \  y - p_k(x) \ \big) ,
\]  where $p_k(t) = \sum_{i=1}^k \tfrac{1}{i!} t^i$ is
the
$k^{\text{th}}$ Taylor polynomial of $e^t - 1$.

\vskip 5pt
\item[(vi).] Assume that $X$ is affine (and as always
non-singular), so that ideal sheaves are identified with
ideals in the coordinate ring $\CC[X]$ of $X$. Given
any non-zero ideal $\fra \subseteq \CC[X]$, set
\[ \fra^{ \{ k \} } \ = \ \big \{ f \in \CC[X] \ \mid
\ Df
\in \fra \ \, \forall \text{ differential operators }
D \text{ on $X$ of order $< k$ } 
\big \}. \]
This determines a  graded family $ \fra^{ \{
\bullet \} } $
 which also reduces to the symbolic  powers $\{
\frq^{(k)}
\}$ when
$\fra =
\frq$ is radical. 

\vskip 5pt
\item[(vii).] Let $\fra_{\bullet} = \{ \fra_k \}$ be a
graded family, and $\frb \subseteq \OO_X$ a fixed
ideal. Then the colon ideals
\[ \frr_k \ = \ \big( \fra_k \ : \ \frb^k \big) \
=_{\text{def}} \ \big \{ f \in \OO_X \mid f \cdot
\frb^k \in \fra_k \, \big \} \] form a graded family.

\vskip 5pt

\end{itemize}
\end{example}

We now construct the asymptotic multiplier ideal
associated to a graded family $\fra_{\bullet}$.
\begin{lemma} \label{Asympt.MI.Lemma} Let
$\fra_{\bullet} =
\{ \fra_k
\}$ be a graded family of ideals, and fix $\ell \in \NN$
plus a rational number $c > 0$. Then for all
positive integers
$p , n
\ge 1$ one has
\begin{equation} \label{gri.inclusions}
\MI{ \tfrac{c}{p} \cdot \fra_{p \ell }} \subseteq
\MI{\tfrac{c}{pn} \cdot \fra_{pn \ell}}.  \notag
\end{equation}
\end{lemma}

\begin{proof} 
Let $\mu : X^{\pr} \lra X$ be a common log resolution
of $\fra_{p\ell}$ and $\fra_{pn\ell}$, with
\[ \fra_{p\ell} \cdot \OO_{X^\pr} \ = \
\OO_{X^\pr}(-F_{p\ell}) \ \ \ , \ \ \fra_{p n \ell}
\cdot
\OO_{X^\pr}
\ = \
\OO_{X^\pr}(-F_{pn\ell}).\]
Condition (\ref{gsi.cond}) implies that
$\fra_{p\ell }^n
\subseteq  \fra_{p n  \ell }$, and hence $- n
F_{p \ell}
\preccurlyeq - F_{p n \ell}$ (i.e. the difference $nF_{p
\ell} - F_{p n \ell}$ is effective). Therefore 
\[
K_{X^{\pr}/X} - [\tfrac{c}{p \ell} \cdot F_{p \ell}]
\ \preccurlyeq \ K_{X^{\pr}/X} - [\tfrac{c}{p
n \ell  }\cdot F_{p n \ell }], 
\]
and the statement follows. 
\end{proof}

 We assert next that  the
collection of multiplier ideals 
\begin{equation} \label{fam.ideals} \Big \{
\MI{\tfrac{c}{p}
\cdot \fra_{p\ell}} \Big \}_{(p > 0)} 
\end{equation} has a unique
maximal element. In fact, the existence of one maximal
element follows from the ascending chain condition on
ideals. On the other hand, if $\MI{\tfrac{c}{p} \cdot
\fra_{p\ell}}$ and $\MI{\tfrac{c}{q} \cdot
\fra_{q\ell}}$ are both maximal,  then thanks to the
Lemma they each coincide with
$\MI{\tfrac{c}{pq} \cdot \fra_{pq\ell}}$.
\begin{definition} Given a graded family of ideals
$\fra_{\bullet} = \{ \fra_k \}$, the \textit{asymptotic
multiplier ideal} at level $\ell$ associated to $c > 0$
and 
$\fra_{\bullet}$, written $\MI{c \cdot \Vert \fra_{\ell}
\Vert}$, is the maxmial element of the collection of
ideals appearing in (\ref{fam.ideals}). In other words, 
\begin{equation}  \label{def.ami.eq} \MI{c \cdot \Vert
\fra_{\ell}
\Vert} \ = \ \MI{\tfrac{c}{p} \cdot   \fra_{p\ell}}
\ \for \text{sufficiently divisible } p \gg 0.  \qed
\end{equation}
\end{definition}

 Assuming as we are that $\fra_k \ne (0)$ for $k \gg
0$, one can show that there is an integer $p_0 =
p_0(\fra_{\bullet},
\ell)$ such that
$\MI{c \cdot
\Vert \fra_{\ell} \Vert} = \MI{\tfrac{c}{p} \cdot
\fra_{p{\ell}}}$ for all $p \ge p_0$. We  use
this fact only to observe that one does not actually
need  the divisibility condition in (\ref{def.ami.eq}).

\begin{remark} Note that  $\MI{c \cdot \Vert
\fra_{\ell} \Vert}$  depends not just on the
particular ideal $\fra_{\ell}$, but on all the ideals
$\fra_{p{\ell}}$ for $p
\gg 0$. The double vertical lines should serve
as a reminder of this point. 
\end{remark}

\begin{example} 
\begin{itemize}
\item[(i).]
If $\fra_k = \fra^k$ is the trivial graded family
consisting of powers of a fixed ideal $\fra$, then
$\MI{c \cdot \Vert \fra_{\ell} \Vert} = \MI{c \cdot
\fra^{\ell}} = \MI{c \ell \cdot \fra}$.
\vskip 5pt
\item[(ii).] When $\frb_k = \frb \big(|kD| \big)$ is
the family of base ideals associated to a big
divisor $D$, then $\MI{c \cdot \Vert \frb_{\ell} \Vert}
= \MI{ c \cdot \Vert  \ell D \Vert}$ is the asymptotic
multiplier ideal  constructed for
instance in \cite{Kawamata} and 
\cite{PAG}. These  ideals have
played an important role in recent work on linear
series.
\vskip 5pt
\item[(iii).] Let $\frq \subseteq \OO_X$ be a radical
ideal. We denote the asymptotic multiplier ideal
at level $\ell$ associated to the symbolic powers
$\frq^{(\bullet)} = \{ \frq^{(k)} \}$ by $\MI{c \cdot
\Vert \frq^{(\ell)}\Vert}$. Thus  $\MI{c \cdot
\Vert \frq^{(\ell)}\Vert } = \MI{\tfrac{c}{p} \cdot
\frq^{(p\ell)}}$ for $p \gg 0$. 
\vskip 5pt
\item[(iv).] Consider the ideals $\fro_k \
\subseteq \CC[x,y]$ constructed in  Example
\ref{Examples.GSI} (v)
 associated to the valuation $v(f) = \ord_t \,
f\big(t, e^t - 1\big)$. Then $\MI{ \Vert \fro_\ell
\Vert} = \CC[x,y]$ for every $\ell$. This can be checked
directly using the observation that each $\fro_k$
contains a polynomial whose divisor is a smooth curve.
From a more sophisticated point of view, the
triviality of the multiplier ideal in question is
implied by Theorem B plus the fact that the colength of
$\fro_k$ in
$\CC[x,y]$ grows linearly rather than quadratically in
$k$. 
\end{itemize}
\end{example}

For our purposes the essential properties of these
multiplier ideals  are given by
\begin{proposition} \label{Asympt.MI.Proposition} Let
$\fra_{\bullet} =
\{
\fra_k
\}$ be a graded family  of ideals on the smooth
variety $X$, and fix
$\ell \ge 1$.  Then:
\begin{itemize}
\item[(i).]  $\fra_{\ell} \subseteq \MI{\Vert
\fra_{\ell}
\Vert}$.

\vskip 5pt

\item[(ii).] For every $m \ge 1$ one has the inclusion
\[\MI{\Vert \fra_{m \ell} \Vert} \ \subseteq \ 
\MI{\Vert \fra_{\ell} \Vert }^{\, m}. \] 
\end{itemize}
\end{proposition}

\begin{proof}  
Since the relative canonical bundle $K_{X^\pr / X}$ is
effective, it follows from the definition
that $\fra \subseteq \MI{\fra}$ for any ideal $\fra
\subseteq \OO_X$. Then using Lemma
\ref{Asympt.MI.Lemma} we find that 
\[ \fra_\ell \ \subseteq \
 \MI{\fra_{\ell}} \ \subseteq \
\MI{\tfrac{1}{p} \cdot \fra_{p \ell}}. \]
Taking $p \gg 0$, this gives (i). For (ii), fix $p \gg
0$ and use the subadditivity relation
(\ref{Subadditivity.additive}) to deduce:
\begin{align*}
\MI{\Vert \fra_{m \ell} \Vert} \ &= \ \MI{ \tfrac{1}{p}
\cdot
\fra_{pm \ell} }  \\
 &= \ \MI{\tfrac{m}{pm} \cdot \fra_{pm \ell} } \\
&\subseteq \   \MI{\tfrac{1}{pm} \cdot \fra_{pm \ell}}^m
\\
 &=   \MI{\Vert \fra_{\ell} \Vert}^ m,
\end{align*}
as asserted.
\end{proof}

\begin{remark} 
Note that it need not be true in general that $\MI{
\fra_{m\ell}} \subseteq \MI{\fra_\ell}^m$. This
explains why it is crucial to pass to the asymptotic
  ideals. 
\end{remark}

\section{Applications}

 Our concrete results  follow from the following
general statement --- which appears as Theorem
B in the Introduction --- concerning the multiplicative
behavior of graded families of ideals:
\begin{theorem} \label{Abstract.Theorem}
Let $\fra_{\bullet} = \{ \fra_k \}$ be a graded family
of ideals on a smooth complex variety $X$, and suppose
that $\frb \subseteq \OO_X$ is an ideal such that
$\MI{\Vert \fra_{\ell} \Vert} \subseteq \frb$ for some
index $\ell \in \NN$. Then
\[ \fra_{ m \ell} \subseteq \frb^m \]
for every integer $m \ge 1$.
\end{theorem}

\noi \textit{Proof.} This is an immediate consequence
of Proposition \ref{Asympt.MI.Proposition}, which
implies that 
\[ \fra_{m \ell} \ \subseteq \ \MI{\Vert \fra_{m \ell}
\Vert} \ \subseteq \ \MI{\Vert  \fra_{\ell} \Vert}^m .
\qed
\]

The first application is to symbolic
powers:\footnote{See \S 3 for a more general
statement.}
\begin{theorem} \label{Symb.Power.Thm} Let $X$ be a
smooth complex variety, and $Z
\subseteq X$  a reduced subscheme all of whose
irreducible components have codimension $\le e$ in $X$.
Put $\frq =
\II_Z$, and fix an integer   $\ell \ge e$. Then
\[
\frq^{(m\ell)} \ \subseteq \big( \frq^{(\ell + 1 - e)}
\big)^m \]
for every $m \ge 1$. In particular, taking $\ell = e$
one has
\[ \frq^{(m e )} \subseteq \frq^m  \fall m \ge 1. \]\
\end{theorem}
\begin{proof}
It suffices by Theorem \ref{Abstract.Theorem} to show
that 
\begin{equation}
\MI{\Vert \frq^{(\ell)}\Vert } \ \subseteq \ \frq^{(\ell
+ 1 - e)}
\tag{*}.
\end{equation}
But membership in the ideal on the right is tested
locally at a general point of each irreducible
component of $Z$. So we can assume after shrinking
$X$ that
$Z$ is smooth and irreducible, of codimension $e$, and
in this case (*) is clear. For then $\frq^{(k)} =
\frq^k$ for all $k$, and $\frq$ is resolved by
taking
$\mu : X^{\pr} =
\Bl_{Z}(X) \lra X$ to be the blow-up of $X$ along
$\frq$.   Writing $E \subseteq X^\pr$ for the
corresponding exceptional divisor, one has 
\[  K_{X^\pr / X}  =
(e-1)E   \ \ \text{and} \ \ 
\frq^{\ell} \cdot \OO_{ X^\pr} = \OO_{X^\pr}(-\ell E).
\]
Consequently
\[\MI{\Vert \frq^{(\ell)}\Vert} \ = \ \mu_*
\OO_{X^\pr}\big ( K_{X^\pr / X} -
\ell E
\big)
\ = \ \mu_*\OO_{X^\pr}\big( - (\ell + 1 - e) E \big) \ =
\frq^{l + 1 - e}, \]
as asserted. \end{proof}

\begin{example} \label{Pts.in.Plane} The first
non-trivial case of Theorem
\ref{Symb.Power.Thm} is the following. Let $T \subseteq
\PP^2$ be a finite set of points, viewed as a reduced
algebraic subset of the plane, and let 
$I
\subseteq
\CC[X, Y, Z]$ be the homogeneous ideal of $T$. If $F
\in \CC[X, Y , Z]$ is a homogeneous polynomial having
multiplicity $\ge 2m$ at every point of $T$, then $F \in
I^m$. (Apply Theorem  \ref{Symb.Power.Thm}  to the
affine cone over $T$ in
$\CC^3$.) In spite of the very
classical nature of this statement  we do not know a
direct elementary proof. 
\end{example}

\begin{remark}
The statement of Theorem \ref{Symb.Power.Thm} can fail
on singular varieties. For example Huneke points out
that counter-examples arise already when  
$Z$ is a line on a quadric cone $X$ in
$\CC^3$.  However Hochster and Huneke \cite{HH} give
some statements valid also on singular ambient spaces. 
\end{remark}

\begin{remark}  
Using familiar arguments, one can deduce from Theorem
\ref{Symb.Power.Thm} that the corresponding statement
holds for excellent regular local rings containing a
field of characteristic zero.  However 
Hochster and Huneke \cite{HH} have shown that in fact
the analogue of (\ref{Symb.Power.Thm}) holds in any 
regular local ring containing a field.  Therefore we
do not dwell on the question  finding the most general
situation in   which the arguments of the present paper
apply. 
\end{remark}

We conclude with a result which renders effective
and extends in certain directions a formulation due to
H\"ubl and Swanson    (\cite{Hub.Swan}, (1.4)) of a
theorem of Izumi
\cite{Iz}:
\begin{theorem} \label{Izumi}
Let $\nu : Y \lra X$ be a proper birational map between
smooth complex varieties. Let
$E \subseteq Y$ be a prime divisor, set \[ \ell \ = \
1 + \text{ord}_E ( K_{Y/X}) \] and for $k \ge 1$ put
$\fro_k = 
\nu_*\OO_Y(-kE)$. Fix an irreducible subvariety  $Z
\subseteq X$  such that $Z \subseteq
\nu(E)$ and denote by $\frp = \II_Z$ the ideal of $Z$.
Then 
\[ \fro_{ \ell m } \ \subseteq \  \frp^m \fall m \ge 1.
\]
\end{theorem}
\begin{remark}
The result discussed in \cite{Hub.Swan} --- which  holds
in  considerably more general  settings, but without
the explicit determination of the coefficient $\ell$ of
$m$ --- deals with the situation in which
$E$ maps to a point. It was in trying to understand
this result that we were led to the statements about
symbolic powers. 
\end{remark}
\begin{proof} [Proof of Theorem \ref{Izumi}]
We can assume without loss of generality that $E$ is
$\nu$-exceptional and that $Z =
\nu(E)$, so that $\nu_* \OO_Y(-E) = \frp$.  Applying 
(\ref{Abstract.Theorem}) to the graded family
$\fro_{\bullet} = \{
\fro_k \}$ (Example \ref{Examples.GSI}(iv)),    it 
suffices to prove that 
$\MI{\Vert
\fro_\ell \Vert} \subseteq \frp$. We   suppose to
this end that we've fixed a large integer $p \gg 0$ such
that the multiplier ideal
$\MI{\Vert \fro_\ell \Vert} = \MI{ \tfrac{1}{p}
\cdot \fro_{p
\ell}}$
 in question is computed on a log resolution
$\mu : X^\pr
\lra X$ of
$\fro_{p
\ell}$ dominating $\nu: Y \lra X$. Then $E$ gives rise
to a prime divisor on $E^\pr$ on $X^\pr$ --- viz.
the proper transform of $E$ --- with
\[  \text{ord}_{E^\pr}(K_{X^\pr / X}) \ = \ 
\text{ord}_{E}(K_{Y/ X}) \ = \ \ell - 1,\] and one
has 
$\fro_{k} =
\mu_* \OO_{X^\pr}(-kE^\pr)$ for every $k \ge 1$.
 
Let $F$ be the effective Cartier divisor on $X^\pr$
defined in the usual way by writing
$\fro_{p \ell} \cdot
\OO_{X^\pr} \  = \
\OO_{X^\pr}\big(-F\big)$. 
Since $\fro_{p \ell} = \mu_* \OO_{X^\pr}\big(- p \ell
E^\pr\big)$, we see that 
$E^\pr$ appears with coefficient $\ge p \ell $ in $F$. 
Consequently
\[ \text{ord}_{E^{\pr}} \Big( K_{X^\pr / X}
- \big[ \tfrac{1}{p}  F \big] 
\Big )
\ \le
\ (\ell -1) -   \ell  \ =  \ -1, \] and therefore
\[ \MI{\Vert \fro_{\ell} \Vert} \ = \ \mu_* \OO_{X^\pr}
\big(  K_{X^\pr / X} -  [ \tfrac{1}{p}  F  ] \big) \
\subseteq
\  \mu_*
\OO_{X^\pr}\big( - E^\pr \big) \ = \ \frp, \]
as required. 
\end{proof}

\section{Generalizations}

 In the preprint \cite{HH}, which appeared shortly
after the first version of the present paper, 
Hochster and Huneke use the theory of tight closure to
extend Theorem \ref{Symb.Power.Thm} to arbitrary
 regular local rings containing a field. They
also observe that it is sufficient to assume that $\frq
\subseteq \OO_X$ is unmixed. In this
section we indicate 
how one applies Theorem
\ref{Abstract.Theorem} to treat unmixed ideals. 

We start by recalling the definition of symbolic
powers in  this more general setting. Assume for
simplicity of exposition that the smooth complex
variety 
$X$ is affine. Given  an ideal
$\frq
\subseteq
\CC[X]$,  fix a primary decomposition
\begin{equation} \frq \ = \ \frq_1
\cap \ldots \cap \frq_h \tag{*} \end{equation}
of $\frq$, and let $Y_i = \Zero{ \sqrt{\frq_i}}$ be the
subvarieties of $X$ corresponding to the associated
primes  $\frp_i = \sqrt{\frq_i}$ of $\frq$. Recall  that
$\frq$ is \textit{unmixed} if  none of the associated
primes $\frp_i$ are embedded (or equivalently if
there are no inclusions among the $Y_i$).  
Then the symbolic powers
$\frq^{(k)} \subseteq \CC[X]$ of $\frq$ are defined as
follows. For each  associated subvariety $Y_i$ of
$\frq$, there is a natural map $ \phi_i : \CC[X] \lra
\OO_{Y_i}X$ from the coordinate ring of $X$ to the local
ring of $X$ along
$Y_i$. We then set
\[  \frq^{(k) } \ = \ \bigcap_{i = 1}^h \ \phi_i^{-1}
\Big (\frq^k
\cdot \OO_{Y_i}X \Big ). \]
In other words, $f \in \frq^{(k)}$ if and only if there
is an element $s \in \CC[X]$, not lying in any of the
associated primes $\frp_i$ of $\frq$, such that $fs \in
\frq^k$. 

Theorem \ref{Symb.Power.Thm} then admits the following
\begin{variantnn} Let $\frq \subseteq
\CC[X]$ be an unmixed ideal, and assume that every
associated subvariety $Y_i$ of $\frq$ has codimension
$\le e $ in
$X$. Then $\frq^{(me)} \subseteq \frq^m$ for all
natural numbers $m \ge 1$. 
\end{variantnn}
\begin{proof} [Sketch of Proof]
The symbolic powers $\frq^{(\bullet)} = \{ \frq^{(k)}\}$
again form a graded family of ideals, so  Theorem
\ref{Abstract.Theorem} will apply as soon as we
establish that $\MI{ \Vert \frq^{(e)} \Vert} \subseteq
\frq$. Referring to the primary decomposition (*), it
is enough to show that 
\begin{equation} \MI{ \Vert \frq^{(e)} \Vert} \subseteq
\frq_i \ \ \text{for each}  \ 1 \le i \le h \tag{**}. 
\end{equation}
For a given index $i$, inclusion in $\frq_i$ is tested
at a generic point of $Y_i$. So having fixed $i$ we are
free to replace $X$ by any open subset meeting $Y_i$.
Therefore, by definition of the symbolic powers, we may
assume after localizing that $\frq^{(k)} =
\frq^{k}$. But in this case $\MI{\Vert \frq^{(e)}
\Vert} = \MI{ \frq^e}$, and  $\MI{\frq^e}
\subseteq
\frq \subseteq \frq_i$ thanks to a variant of a theorem
of  Skoda (cf.
\cite{MIAG}). 
\end{proof}

\begin{remark}  While there are certain similarities of
spirit between the arguments appearing here and those
of Hochster and Huneke \cite{HH}  --- e.g. both
involve asymptotic constructions, and reduce to the
situation in which
$\frq^{(k)} =
\frq^k$ --- the precise connections between the
two points of view remain quite mysterious.  In the
hopes of understanding these connections more clearly,
it is interesting to observe that the properties of
multiplier ideals used here can
be ``axiomatized" as follows. Given a graded family
$\frq_{\bullet} = \{
\frq_k \}$ what is required for the application to
symbolic powers is the existence of ideals 
$\mathcal{H}(\Vert\frq_m\Vert) \subseteq \OO_X$ 
satisfying  the following properties:
\begin{enumerate}
\item[(i).]  $\mathcal{H}(\Vert\frq_m\Vert)$ is a sheaf
on $X$, i.e. it commutes with localization, and when
$\frq_k = \fra^k$ is the trivial family consisting of
powers of a fixed ideal $\fra$, then  Skoda's theorem
\[  \mathcal{H}(\Vert\fra^n\Vert) \subseteq \fra \]
holds;\footnote{One also could ask for more precise
statements involving the codimensions of associated
primes of $\fra$.} 
\item[(ii).] $\frq_m
\subseteq \mathcal{H}(\Vert\frq_m\Vert) \fall m$;
\item[(iii).] One has the subadditivity relation:
\[\mathcal{H}(\Vert\frq_{\ell m}\Vert) \subseteq
\mathcal{H}(\Vert\frq_m\Vert)^\ell.\]
\end{enumerate}
In our setting, the required ideals are of course given
by the asymptotic multiplier ideals $\MI{\Vert \frq_m
\Vert}$.  However the existence of such ideals
$\mathcal{H}$ is a purely algebraic question, and it
would be very interesting to give a 
construction e.g. using ideas from tight closure. The
hope here is that such a construction might serve
as a Rosetta stone to help in deciphering the
connections between the methods of the present note
and the theory of tight closure.
\qed
\end{remark}


\begin{thebibliography}{99}


\bibitem{Demailly} J.-P.\ Demailly, {\em $L^2$
vanishing theorems for positive line bundles and
adjunction theory}, Lecture Notes of the CIME Session
``Transcendental methods in Algebraic Geometry'',
Cetraro, Italy, July 1994, Ed.\ F.~Catanese,
C.~Ciliberto, Lecture  Notes in Math., Vol.~1646, 1--97.

\bibitem{Demailly.Bourbaki} J.-P.\ Demailly, {\em
M\'ethodes $L^2$ et  r\'esultats effectifs en
g\'eom\'etrie alg\'ebrique}, S\'eminaire Bourbaki, 
expos\'e n${}^\circ\,$852, novembre 1998, 38pp.


\bibitem{DEL} J.-P. Demailly, L. Ein and R. Lazarsfeld,
\textit{A subadditivity property of multiplier ideals},
to appear in the special issue of the  Mich. Math. J.
in  honor of W. Fulton.


\bibitem{Ein}  L. Ein, {\em Multiplier ideals,
vanishing theorems and applications}, in
\textit{Algebraic Geometry -- Santa Cruz 1995},  Proc.
Symp. Pure Math. \textbf{62} (1997), pp. 203 -- 219. 

\bibitem{Eisenbud} D. Eisenbud, {\em
Commutative Algebra, with a view towards algebraic
geometry},  Graduate Texts in Math.  {\bf 150}, 
Springer-Verlag, New York, 1995.


\bibitem{HH} M. Hochster and C. Huneke, {\em Comparison
of symbolic and ordinary powers of ideals}, preprint.


\bibitem{Hub.Swan} R. H\"ubl and I. Swanson,
{\em Discrete valuations cenetered on local domains},
preprint (1998).

\bibitem{Huneke} C. Huneke, {\em Uniform bounds in
Noetherian rings}, Inv. Math. \textbf{107} (1992), pp.
202 - 223.


\bibitem{Iz} S. Izumi, \textit{A measure of integrity
for local analytic algebras}, Publ. RIMS Kyoto Univ.
\textbf{21} (1985), pp. 719--735. 

\bibitem{Kawamata} Y.\ Kawamata, {\em Deformations of
canonical singularities}, Journal of the Amer.\ Math.\
Soc.\ {\bf 12} (1999), 519--527.


\bibitem{PAG} R. Lazarsfeld, \textit{Positivity in
Algebraic Geometry}, book in preparation.

\bibitem{MIAG} R. Lazarsfeld, \textit{Multiplier ideals
for algebraic geometers}, informal lecture notes
available at
\texttt{http://www.math.lsa.umich.edu/\~{}rlaz}

\bibitem{Lipman} J. Lipman, \textit{Adjoints of ideals
in regular local rings}, Math. Res. Lett. \textbf{1}
(1994),  pp. 739 -- 755.

\bibitem{Siu} Y.-T.\ Siu, {\em Invariance of
plurigenera}, Invent.\ Math.\ {\bf 134} (1998),  661--673.

\bibitem{Swanson} I. Swanson, {\em
Linear equivalence of topologies}, preprint (1997).
\end{thebibliography}
\end{document}